\documentclass[journal]{IEEEtran}
\ifCLASSINFOpdf
\else
\fi
\usepackage{amsmath}
\usepackage{makecell}
\usepackage{rotating}
\usepackage{ctable}
\usepackage[graphicx]{realboxes}
\usepackage{amsthm}
\usepackage{amsmath}
\usepackage{amssymb}
\usepackage{ulem}
\newtheorem{thm}{Theorem}

\newtheorem{dfntn}{Definition}[section]

\newcommand{\pdf}[1]{f_{#1}}
\newcommand{\cep}[2]{\hat{c}_{#1}[#2]}
\newcommand{\cumu}[2]{\kappa_{#1}[#2]}

\newcommand{\cf}[1]{\Phi_{#1}(\mu)}
\newcommand{\prob}[1]{\mathbb{P}(#1)}

\newcommand{\cbg}[1]{\textcolor{red}{(BG: #1)}}

\renewcommand{\sout}[1]{}

\begin{document}
\title{Cepstral Analysis of Random Variables: Muculants}

\author{Christian~Knoll,~\IEEEmembership{Student Member,~IEEE},~Bernhard~C.~Geiger,~\IEEEmembership{Member,~IEEE}, and~Gernot~Kubin,~\IEEEmembership{Member,~IEEE}
\thanks{The authors are with the Signal Processing and Speech Communication Laboratory, Graz University of Technology, 8010 Graz, Austria. (e-mail: christian.knoll.c@ieee.org, geiger@ieee.org, g.kubin@ieee.org).}}%
\maketitle
\begin{abstract} 
An alternative parametric description for discrete random variables, called \emph{muculants}, is proposed. In contrast to cumulants, muculants are based on the Fourier series expansion, rather than on the Taylor series expansion, of the logarithm of the characteristic function. We utilize results from cepstral theory to derive elementary properties of muculants, some of which demonstrate behavior superior to those of cumulants. For example, muculants and cumulants are both additive. While the existence of cumulants is linked to how often the characteristic function is differentiable, all muculants exist if the characteristic function satisfies a Paley-Wiener condition. Moreover, the muculant sequence and, if the random variable has finite expectation, the reconstruction of the characteristic function from its muculants converge. We furthermore develop a connection between muculants and cumulants and present the muculants of selected discrete random variables. Specifically, it is shown that the Poisson distribution is the only distribution where only the first two muculants are nonzero.
\end{abstract}


%
\IEEEpeerreviewmaketitle

\section*{Introduction}
\IEEEPARstart{C}{umulants} are parametric descriptors of random variables (RVs) and are commonly used to analyze non-Gaussian processes~\cite{mendelHOS} or the effects of nonlinear systems~\cite{petropulu}. Unlike moments, cumulants are orthogonal descriptors and satisfy a homomorphism property~\cite{moments-cumulants}.

As attractive as these properties are, there are several drawbacks: 
First, since cumulants are the Taylor series coefficients of the logarithm of the distribution's characteristic function~\cite{mendelHOS}, they constitute a complete description only for infinitely differentiable  characteristic functions.
Second, there are no general results regarding the behavior of the sequence of cumulants;
the sequence may even diverge (this problem, though, can be mitigated by the definition of q-moments and q-cumulants~\cite{qMoments, qCumulants}). Hence, a reconstruction of a distribution function in terms of its cumulants~\cite{cramer1957,petrov1972,edgeworthExpansion} may not converge.
Third, the Marcinkiewicz theorem~\cite{char_func_luk} states that any RV either has infinitely many cumulants or up to second order only;
hence, every cumulant-based approximation is either Gaussian or does not correspond to a RV at all. It follows that cumulants are not well suited for hypothesis testing, except in the important case of testing Gaussianity.
Fourth, if the RV takes values in the set of integers, then the characteristic function is periodic, and the Taylor series expansion fails to remain the most natural approach.

In this paper, motivated by the latter shortcoming,
we replace the Taylor series expansion by a Fourier transform. While in general the resulting description is a functional, for integer RVs, the Fourier transform degenerates to a Fourier series expansion. The resulting Fourier coefficients -- henceforth called \emph{muculants} -- are a parametric (i.e., finite or countable) description, retaining several properties of cumulants while removing several of their shortcomings. For example, muculants are orthogonal and additive, but truncating the series leads to a bounded approximation error that converges to zero with increasing order.
Also the existence of muculants is less problematic than the one of cumulants, as the former is ensured if the characteristic function satisfies a Paley-Wiener condition. Finally, even though the Poisson distribution is the only distribution for which only the first two muculants are nonzero (see Section~\ref{sec:Experiments} for the muculants of selected discrete distributions), there exist distributions with more than two (but still only finitely many) nonzero muculants. 

The sequence of operations -- Fourier transform, (complex) logarithm, and inverse Fourier transform -- is also essential in cepstral analysis, originally introduced to investigate the influence of echo~\cite{Bogert_1963} and to represent nonlinear systems~\cite{oppenheim-phd}. Today the cepstrum is widely used in speech processing~\cite{oppenheim_quefrency}. Our analysis of RVs is thus deeply rooted in signal processing. Existence and properties of muculants, stated in Section~\ref{sec:muculants}, are based on properties of the cepstrum, or more generally, properties of the Fourier transform. Moreover, a connection between muculants and cumulants, presented in Section~\ref{sec:Cumulants} also finds counterparts in cepstral analysis~\cite{schroeder-cepstrum-predictor, khare-moments}.

\sout{In this sense, our work follows the spirit of Widrow, who also applied methods from signal processing to the analysis of RVs~\cite{widrow_phd,widrow_quantization}.\footnote{Specifically, the nonlinear quantization operation is a linear sampling operation of the probability density function of the quantizer input signal. This allows for a quantization theorem reminiscent of Shannon's sampling theorem, stating conditions under which the distribution of the quantizer output can be used to reconstruct the distribution of its input.}}

About terminology: The name \emph{cepstrum} is derived from reversing the first syllable of the word \emph{spectrum}: While the cepstrum is situated in the original (e.g., time) domain, this terminology was introduced to emphasize that this sequence of operations can provide fundamentally different insights~\cite{Bogert_1963}. Following this approach, we call our parametric descriptors \emph{muculants}, a reversal of the first syllable of \emph{cumulants}.

\section{Preliminaries} \label{sec:preliminaries}
Let $\textbf{X}$ be a real-valued RV and let $\cf{\mathbf{X}}$ denote its characteristic function
\begin{equation}\label{eq:def:cf}
 \cf{\mathbf{X}} := \mathbb{E}\left(e^{\jmath \mu \textbf{X}}\right),
\end{equation}
 where $\mu \in \mathbb{R}$.
The characteristic function always exists since it is an expectation of a bounded function (see~\cite{char_func_luk,luk-survey-cf,applications-luk} for the theory, properties, and applications of characteristic functions). Two RVs are identically distributed iff (if and only if) their characteristic functions are equal. Moreover, every characteristic function satisfies the following properties~\cite[p.~13]{bisgaard2000}: (i) $\cf{\mathbf{X}}$ is uniformly continuous everywhere, (ii) $\cf{\mathbf{X}}|_{\mu=0}=1$, (iii) $|\cf{\mathbf{X}}| \leq 1$, and (iv)  $\cf{\mathbf{X}}$ is Hermitian. Finally, if $\mathbf{X}$ and $\mathbf{Y}$ are \emph{independent}, then the characteristic function of $\mathbf{Z}=\mathbf{X}+\mathbf{Y}$ is $\cf{\mathbf{{Z}}} = \cf{\mathbf{X}} \cdot \cf{\mathbf{Y}}$.

Let $\textbf{X}$ be a discrete RV taking values from the set of integers $\mathbb{Z}$. It can be described by its probability mass function (PMF) $\pdf{\mathbf{X}}{:}\ \mathbb{Z}\to [0,1]$, where
\begin{equation}
  \forall \xi\in\mathbb{Z}: \quad\pdf{\mathbf{X}}[\xi] := \prob{\textbf{X}=\xi}.
\end{equation}
In this case,~\eqref{eq:def:cf} equates to $\cf{\mathbf{X}} = \sum_{\xi=-\infty}^{\infty} \, \pdf{\mathbf{X}}[\xi] e^{\jmath \mu \xi}$,
i.e., $\cf{\mathbf{X}}$ is the inverse Fourier transform of $\pdf{\mathbf{X}}$ and periodic.

We call a PMF \emph{causal} if $\pdf{\mathbf{X}}[\xi]=0$ for $\xi<0$. We call a PMF \emph{minimum-phase} if it is causal and its $z$-transform
\begin{equation}
 \Psi_{\mathbf{X}}(z) := \sum_{\xi=0}^\infty \pdf{\mathbf{X}}[\xi] z^{-\xi}
\end{equation}
has all its zeros inside the unit circle.\footnote{Possible poles are inside the unit circle by construction, since every PMF is absolutely summable.} The characteristic function is related to this $z$-transform via $\cf{\mathbf{X}} = \Psi_{\mathbf{X}}\left(e^{\jmath \mu}\right)$.

For a complex function $z(t)$ without zeros on the unit circle, the complex logarithm is uniquely defined as
\begin{equation}
	\log z(t) = \ln|z(t)| +\jmath \cdot \textrm{arg}^*(z(t)),
	\label{eq:complex-log}
\end{equation}
where $\ln$ is the natural logarithm and $\textrm{arg}^*(z(t)) = \textrm{Arg}(z(t)) + 2k\pi$ is continuous for $k\in\mathbb{Z}$~\cite{tribolet-1977};
i.e., one formally 
represents $\log(z(t))$ as a single-valued analytic function on a Riemann surface~\cite{brown-complex}. The computation of such a continuous phase function is essential for the estimation of the complex cepstrum in practical cases ~\cite{oppenheim2010discrete, tribolet-1977, mc-gowan, karam}.
In contrast to the principal value of the complex logarithm (cf.~\cite[p.~1009]{oppenheim2010discrete}), the complex logarithm defined in~\eqref{eq:complex-log} satisfies, for complex functions $w(t)$ and $z(t)$,
\begin{equation}\label{eq:logproperties}
 \log(w(t)z(t)) = \log(w(t)) + \log(z(t)).
\end{equation}


Finally, the \emph{cumulants} $\{\cumu{\mathbf{X}}{n}\}_{n\in\mathbb{Z}}$ are the coefficients of the Taylor series expansion
of $\log \cf{\mathbf{X}}$, i.e.,
\begin{equation}\label{eq:cumulant:def}
  \log \cf{\mathbf{X}} = \sum_{n=1}^{\infty} \cumu{\mathbf{X}}{n} \frac{(\jmath \mu)^n}{n!}
\end{equation}
provided the sum on the r.h.s.\ exists. Specifically, if $\mathbb{E}(\mathbf{X}^n)<\infty$, then $\cf{\mathbf{X}}$ is $n$ times continuously differentiable~\cite[p.~48]{bisgaard2000}, and we obtain the $n$-th cumulant as
\begin{equation}\label{eq:cumulant:derivative}
 \cumu{\mathbf{X}}{n} = \frac{d^n \log\cf{\mathbf{X}}}{\jmath^nd\mu^n} \Bigm|_{\mu = 0}.
\end{equation}
If $\mathbf{X}$ and $\mathbf{Y}$ are \emph{independent}, then $\cumu{\mathbf{X}+\mathbf{Y}}{n}=\cumu{\mathbf{X}}{n}+\cumu{\mathbf{Y}}{n}$, i.e., cumulants are additive. For elementary results on cumulants the reader is referred to~\cite{mattner-cumulants} and~\cite{rota-shen}.

\section{Muculants: Definition and Properties}\label{sec:muculants}
The definition of the muculants follows the definitions of the cepstrum~\cite[Ch.~13]{oppenheim2010discrete}, with the main difference that the roles of Fourier transform and inverse Fourier transform are reversed. While some properties of the muculants are directly transferred from cepstral analysis, several results are based on the fact that we operate on a probability space.

\begin{dfntn}[Complex Muculants]\label{def:complex}
The \textbf{complex muculants} $\{\cep{\mathbf{X}}{n}\}_{n\in\mathbb{Z}}$ are the coefficients of the Fourier series expansion (if it exists) of $\log \cf{\mathbf{X}}$, i.e., 

\begin{equation}
	\log \cf{\mathbf{X}} = \sum_{n=-\infty}^{\infty} \, \cep{\mathbf{X}}{n}e^{\jmath  \mu n}.
\label{eq:complex-cepstrum-sum}
\end{equation}
\end{dfntn}
Note that~\eqref{eq:complex-cepstrum-sum} is a one-to-one correspondence only if $\log$ is uniquely defined, as in~\eqref{eq:complex-log}. If $\cf{\mathbf{X}}$ has zeros on the unit circle, phase jumps by $\pi$ may be introduced and one cannot determine an unambiguous phase curve~\cite{tribolet-1977}. Since~\eqref{eq:complex-cepstrum-sum} must be one-to-one to establish the desirable property of additivity, we consider only RVs for which $\cf{\mathbf{X}}$ has no zeros on the unit circle. We note in passing that these phase ambiguities are unproblematic for \textbf{power muculants}, i.e., the Fourier coefficients of $\ln |\cf{\mathbf{X}}|^2$, which are complete parametric descriptors of minimum-phase PMFs (see~\cite{muculants-arxiv},~\cite[p.~215]{papoulis_fourier}).



With $\log z(t)$ defined as in~\eqref{eq:complex-log}, the complex muculants are:
	\begin{equation}
		\hat{c}_{\mathbf{X}}[n] = \frac{1}{2\pi} \int^{\pi}_{-\pi} \log \cf{\mathbf{X}} e^{-\jmath  \mu n}  \; \mathrm{d} \mu.
		\label{eq:complex-cepstrum-integral}
	\end{equation}
This integral is well-defined iff the integrand is absolutely integrable, i.e., iff $\int_{-\pi}^{\pi} |\log\cf{\mathbf{X}}|d\mu < \infty$. Using~\eqref{eq:complex-log} and applying the triangle inequality, we get
\begin{equation}
	|\!\ln|\Phi_{\mathbf{X}}(\mu)\!|| \!\! \leq  \!\! |\textrm{log}(\Phi_{\mathbf{X}}(\mu)\!)|\!\! \leq \!\!|\!\ln|\Phi_{\mathbf{X}}(\mu)\!||\!\! + \!\!|\textrm{arg}^*\!(\Phi_{\mathbf{X}}(\mu)\!)|.
\label{eq:triangleInequ}
\end{equation}
The phase of $\cf{\mathbf{X}}$ is continuous, thus its integral over the compact set $[-\pi, \pi]$ is finite. Since moreover $|\cf{\mathbf{X}}| \leq 1$, the integral in~\eqref{eq:complex-cepstrum-integral} is well-defined iff
\begin{equation}
	\frac{1}{2\pi} \int_{-\pi}^{\pi} \ln|\Phi_{\mathbf{X}}(\mu)|d\mu > -\infty.
\label{eq:paleyWiener4}
\end{equation}
Note that this condition rules out the existence of muculants for RVs whose characteristic function vanishes on an interval.

\sout{\cbg{We could shorten/remove that paragraph, or only talk about causality (and not on minimum-phase).}} The condition in~\eqref{eq:paleyWiener4} reminds of the Paley-Wiener condition~\cite[p.~423]{papoulis_probability}. Loosely speaking and translated to the language of probabilities, this condition states that iff~\eqref{eq:paleyWiener4} holds, then to $\ln|\cf{\mathbf{X}}|$ a \emph{unique} phase can be associated which guarantees that the corresponding PMF is minimum-phase. Causality of $\pdf{\mathbf{X}}$ is thus a sufficient (but not necessary) condition for~\eqref{eq:paleyWiener4}. 

	
The following theorem, proved in Appendix~\ref{app:proof}, summarizes the main properties of complex muculants:
\begin{thm}[Properties of Complex Muculants]\label{thm:properties}
 Let $\mathbf{X}$ be an RV with PMF $\pdf{\mathbf{X}}$ supported on $\mathbb{Z}$ and let $\{\cep{\mathbf{X}}{n}\}_{n\in\mathbb{Z}}$ be the complex muculants defined in Definition~\ref{def:complex}. Then, the following properties hold:
 \begin{enumerate}
  \item $\cep{\mathbf{X}}{n} \in\mathbb{R}$. \label{prop:real}
  \item $\cep{\mathbf{X}}{0}\le 0$.	\label{prop:negative}
  \item If $\pdf{\mathbf{X}}(\xi)=\pdf{\mathbf{X}}(-\xi)$, then $\cep{\mathbf{X}}{n}=\cep{\mathbf{X}}{-n}$. \label{prop:even}
  \item If $\mathbb{E}\left(\textbf{X}\right) < \infty$, then $\sum_n \cep{\mathbf{X}}{n} =0$ and the series in~\eqref{eq:complex-cepstrum-sum} converges pointwise. \label{prop:pointwise}
  \item $\lim\limits_{n \to \pm\infty}{\cep{\mathbf{X}}{n}} = 0$. If $\mathbb{E}\left(\textbf{X}\right) < \infty$, then $\cep{\mathbf{X}}{n} = \mathcal{O}(1/n)$. \label{prop:tozero}
  \item If $\mathbf{X}$ and $\mathbf{Y}$ are independent, then $\cep{\mathbf{X}+\mathbf{Y}}{n}=\cep{\mathbf{X}}{n}+\cep{\mathbf{Y}}{n}$. \label{prop:additive}
 \end{enumerate}
\end{thm}

%


Note that finite expectation is required in properties~\ref{prop:pointwise} and~\ref{prop:tozero}, because one can construct characteristic functions which are nowhere differentiable. E.g., the Weierstrass function $\sum_{n=1}^{\infty} \, \left(\frac{1}{2}\right)^n \cos(3^n \mu)$ is nowhere differentiable, but continuous everywhere and a valid characteristic function~\cite[p.~47]{bisgaard2000}.

That muculants of an RV with finite expectation must sum to zero makes truncating the Fourier series as problematic as truncating the cumulant expansion, i.e., the truncated series need not correspond to a valid PMF. However, by Parseval's theorem and the fact that $\cep{\mathbf{X}}{n} = \mathcal{O}(1/n)$, the squared error between the true and the approximated characteristic function stays bounded. Hence, muculants may behave better than cumulants when used in functionals of distributions, such as entropy or informational divergence. \sout{It would be interesting to investigate, e.g., if an approximation of negentropy \cbg{negentropy = continuous distribution?} based on muculants is possible and has better convergence properties than one based on cumulants~\cite{negentropy}.} Moreover, while there exists no distribution with finitely many (but more than two) nonzero cumulants, distributions with, e.g., only three nonzero complex muculants exist.

Finally, property~\ref{prop:additive} states that muculants are, just as cumulants, additive descriptors. This retains the benefits of cumulants while eliminating some of their drawbacks particularly problematic with discrete RVs.

\section{Linking Cumulants and Muculants}\label{sec:Cumulants}
The $z$-transform points at a close connection between cumulants and the cepstrum~\cite{schroeder-cepstrum-predictor, khare-moments}, thus suggesting a connection between cumulants and muculants. Suppose that $\log\cf{\mathbf{X}}$ is continuous, that its first $(k-1)$ derivatives are continuous, and that its $k$-th and $(k+1)$-th derivatives are piecewise continuous. We then obtain with~\cite[Th~3.22]{maths-sp11} that
\begin{equation}
 \frac{d^k}{d\mu^k} \log\cf{\mathbf{X}} = \sum_{n = -\infty}^{\infty} (jn)^k \cep{\mathbf{X}}{n} e^{\jmath \mu n}.
\end{equation}
Evaluating the l.h.s.\ at $\mu=0$ yields the $k$-th cumulant $\jmath^k\cumu{\mathbf{X}}{n}$ (if it exists, cf.~\eqref{eq:cumulant:derivative}); evaluating the r.h.s.\ at $\mu=0$ then yields
\begin{equation}
 \cumu{\mathbf{X}}{n}= \sum_{n = -\infty}^{\infty} n^k \cep{\mathbf{X}}{n}, \label{eq:cepcum}
\end{equation}
where, abusing terminology by ignoring the fact that the muculants need not be non-negative, we call the r.h.s.\ the \emph{$k$-th non-central moment of the complex muculants}.

In~\cite{khare-moments}, the moments of the cepstrum are connected to the moments of the original sequence, yielding a recursive formula to compute the former from the latter. An equation similar to~\eqref{eq:cepcum}, called Good's formula~\cite{good}, expresses the cumulants in terms of the moments of a random vector.

\section{Muculants of Selected RVs}\label{sec:Experiments}
In this section, we present the complex muculants for a selection of discrete distributions, summarized in Table~\ref{tab:xyzs}. 

\begin{table*}
  \caption{Muculants $\cep{\mathbf{X}}{n}$ of selected distributions. We also present the cumulants $\cumu{\mathbf{X}}{n}$ and the characteristic functions $\cf{\mathbf{X}}$.}
 \label{tab:xyzs}
\Rotatebox{0}{
 \centering.

\begin{tabular}{|l|l|l|l|l|}
\hline
Distribution 
& $\pdf{\mathbf{X}}[\xi]$
& $\cf{\mathbf{X}}$ 
& $\cumu{\mathbf{X}}{n}$ 
& $\cep{\mathbf{X}}{n}$\\
\hline
Poisson
& $\begin{cases} e^{-\lambda}\frac{\lambda^\xi}{\xi!},& \xi \ge 0\\0,& \text{else}\end{cases}$
& $e^{\lambda(e^{\jmath \mu}-1)}$ 
& $\cumu{\mathbf{X}}{n}=\lambda$ 
& $
\begin{cases}
 -\lambda,& \text{ if } n=0\\
 \lambda,& \text{ if } n=1\\
 0, & \text{ else}
\end{cases}$\\
\hline
Degenerate 
& $\delta[\xi-M]$
& $e^{\jmath \mu M}$ 
& $\cumu{\mathbf{X}}{n}=
\begin{cases} 
 M,&\text{ if } k=1\\
 0,&\text{ else} 
\end{cases}$
& $\begin{cases}0, & \text{if }n=0\\\frac{M}{n}(-1)^{n+1}, & \text{else}\end{cases}$\\
\hline
\makecell[l] {Bernoulli\\ ($p<0.5$)}
& $\begin{cases}
 1-p,& \text{ if } \xi=0\\
 p,& \text{ if } \xi = 1\\
 0,&\text{ else }
\end{cases}$ 
& $1-p+pe^{\jmath \mu}$
& $\cumu{\mathbf{X}}{n+1} = p(1-p) \frac{d\cumu{\mathbf{X}}{n}}{dp}$
&  $
\begin{cases}
 \log(1-p),& \text{ if } n=0\\
 \frac{(-1)^{n+1}}{n}\left(\frac{p}{1-p}\right)^n,& \text{ if } n> 0\\
 0, & \text{ else}
\end{cases}$\\
\hline
\makecell[l] {Bernoulli\\ ($p>0.5$)}
&$\begin{cases}
 1-p,& \text{ if } \xi=0\\
 p,& \text{ if } \xi = 1\\
 0,&\text{ else }
\end{cases}$ 
& $1-p+pe^{\jmath \mu}$
& $\cumu{\mathbf{X}}{n+1} = p(1-p) \frac{d\cumu{\mathbf{X}}{n}}{dp}$
&  $
\begin{cases}
 \log p,& \text{ if } n=0\\
 \frac{(-1)^{n+1}}{n}\left(1+\left(\frac{1-p}{p}\right)^n\right),& \text{ if } n< 0\\
 \frac{(-1)^{n+1}}{n}, & \text{ if } n> 0
\end{cases}$
\\
\hline
Geometric
& $\begin{cases}(1-p)p^{\xi},& \xi \ge 0\\0,&\text{else}\end{cases}$
& $\frac{1-p}{1-pe^{\jmath \mu}}$
&\makecell[l]{ $\cumu{\mathbf{X}}{n+1} = \rho(1+\rho) \frac{d\rho}{d\cumu{\mathbf{X}}{n}}$ \\ $\rho=(1-p)/p$}
& $
\begin{cases}
 \log(1-p),& \text{ if } n=0\\
 \frac{p^n}{n},& \text{ if } n> 0\\
 0, & \text{ else}
\end{cases}$\\
\hline
\makecell[l]{Negative\\ Binomial}
& $\begin{cases}\binom{\xi+N-1}{\xi}(1-p)^N p^{\xi},& \xi \ge 0\\0,&\text{else}\end{cases}$
& $\left(\frac{1-p}{1-pe^{\jmath \mu}}\right)^N$
&\makecell[l] {$N$ times the cumulant \\ of the geometric distribution}
&\makecell[l] {$N$ times the muculant \\ of the geometric distribution} \\ \hline
Binomial
& $\begin{cases}\binom{N}{\xi}p^{\xi}(1-p)^{(N-\xi)},& \xi \ge 0\\0,&\text{else}\end{cases}$
& $\left(1-p+pe^{\jmath \mu}\right)^N$
&\makecell[l] {$N$ times the cumulant \\ of the Bernoulli distribution}
& \makecell[l] {$N$ times the muculant \\ of the Bernoulli distribution}
\\
\hline

\end{tabular}
}
\end{table*}

\subsection{Poisson Distribution}\label{sec:poisson}
With $\lambda>0$, the characteristic function of the Poisson distribution is $\cf{\mathbf{X}} = e^{\lambda(e^{\jmath \mu}-1)}$. All cumulants exist and equal $\lambda$. A fortiori, $\mathbb{E}(\mathbf{X})=\cumu{\mathbf{X}}{1}=\lambda$, and with~\eqref{eq:complex-cepstrum-sum} we can write
\begin{equation}
\textrm{log}(e^{\lambda(e^{\jmath \mu}-1)}) = \lambda e^{\jmath \mu} - \lambda = \sum_{n= -\infty}^{\infty} \hat{c}_{\mathbf{X}}[n] \cdot e^{\jmath \mu n}.
\label{eq:poisson-cc}
\end{equation}
Equating the coefficients yields $\cep{\mathbf{X}}{0}=-\lambda$, $\cep{\mathbf{X}}{1}=\lambda$, and $\cep{\mathbf{X}}{n}=0$ for $n\neq 0,1$. With property~\ref{prop:pointwise} of Theorem~\ref{thm:properties}, 
the Poisson distribution is thus the only distribution with all but the first two muculants being zero.

\subsection{Degenerate Distribution}
The PMF for  $\mathbf{X}\equiv M$ is the Kronecker delta, i.e., $\pdf{\mathbf{X}}[\xi]=\delta[\xi-M]$ and $\cf{\mathbf{X}}=e^{\jmath \mu M}$. While $\cep{\mathbf{X}}{0}=0$, for $n\neq 0$ the complex muculants are given as 
\begin{align}
\cep{\mathbf{X}}{n} &= \frac{1}{2\pi}\int_{-\pi}^\pi \jmath \mu M e^{-\jmath \mu n} d\mu = \frac{M}{n}(-1)^{n+1}.
\label{eq:cc-delta1}
\end{align}
Since in this case $\log\cf{\mathbf{X}}=\jmath \mu M$ is not periodic, one cannot expect~\eqref{eq:cepcum} to hold: Indeed, while $\kappa_3=0$, the third non-central moment of the complex muculants diverges.

\subsection{Minimum-Phase Distribution}
Suppose that $\pdf{\mathbf{X}}$ is minimum-phase. Its $z$-transform
\begin{equation}
  \Psi_{\mathbf{X}}(z) = A\cdot\frac{\prod_{k=1}^{\infty}(1-o_k z^{-1})}{\prod_{k=1}^{\infty}(1-p_k z^{-1})}
  \label{eq:poles-zeros}
\end{equation}
has all poles $p_k$ and zeros $o_k$ inside the unit circle, i.e., $|o_k|<1$ and $|p_k|<1$ for all $k$. This applies, for example, to the geometric distribution, the Bernoulli distribution and the binomial distribution with $p<0.5$, and the negative binomial distribution. Exploiting $\cf{\mathbf{X}}=\Psi_{\mathbf{X}}(e^{\jmath \mu})$,~\eqref{eq:logproperties}, and the Mercator series, which for  $|z|\le 1$ and $z\neq -1$ reads
\begin{equation}\label{eq:mercator}
	 \log(1+z) = \sum_{m=1}^\infty \frac{(-1)^{m+1}}{m} z^m,
	\end{equation}
the complex muculants are obtained by, cf~\cite[p.~1011]{oppenheim2010discrete}
\begin{equation}
 \log \cf{\mathbf{X}} = \underbrace{\log A}_{=\cep{\mathbf{X}}{0}} + \sum_{n=1}^\infty \underbrace{\sum_{k=1}^\infty \frac{p_k^n-o_k^n}{n}}_{=\cep{\mathbf{X}}{n}} e^{-\jmath \mu n}.
\end{equation}
Specifically, for minimum-phase PMFs, we obtain $\cep{\mathbf{X}}{n}=0$ for $n<0$. Note further that~\eqref{eq:mercator} can be used to derive a recursive relation among the complex muculants which, at least under special conditions, admits computing the complex muculants from the PMF without requiring a Fourier transform or a complex logarithm~\cite[p.~1022]{oppenheim2010discrete}.

\section{Discussion and Future Work}\label{sec:discussion}
We have argued that truncating the muculant series~\eqref{eq:complex-cepstrum-sum} results in a bounded error that decreases with increasing number of summands, a property that does not hold for the cumulant series~\eqref{eq:cumulant:def}. Although this is an advantage of muculants over cumulants, the question under which circumstances a truncated muculant series represents a PMF remains open. One may investigate, for example, the approximation of a distribution by one with a given number of nonzero muculants. A second line of research could investigate expressions for functionals of discrete distributions (such as entropy and informational divergence) based on muculants, thus complementing cumulant-based expressions for continuous distributions, cf.~\cite{negentropy}.

That the Poisson distribution has only two nonzero muculants (cf.~Section~\ref{sec:poisson}) makes the presented theory attractive for hypothesis testing. Future work shall thus develop numerical methods to estimate muculants from data, together with estimator variances and confidence bounds. Based on that, a test whether a distribution is Poisson or not shall be formalized and compared to existing hypothesis tests.

We presented a condition for the existence of muculants in~\eqref{eq:paleyWiener4}; however, it is not clear whether there exist distributions violating this condition. The search for a concrete example, or preferably, for a more general statement about existence is within the scope of future work. Similar investigations shall be devoted to the convergence of~\eqref{eq:complex-cepstrum-sum}, which for the moment is guaranteed only for distributions with finite expectation (cf.~property~\ref{prop:pointwise} of Theorem~\ref{thm:properties}).

Finally, the theory of muculants shall be extended to continuous RVs, even though this requires a \emph{muculant function} rather than a muculant sequence. In this context, the connection between the muculant function and/or cumulants of a continuous RV and the muculants of a discrete RV obtained by uniformly quantizing the former might be of interest. A fundamental step towards these results lies in the fact that quantization periodically repeats the characteristic function~\cite{widrow_quantization}.

\appendices
\section{Proof of Theorem~\ref{thm:properties}}\label{app:proof}
 \emph{1)} The Fourier transform of a Hermitian function is real-valued~\cite[p.~83]{oppenheim2010discrete}; $\cf{\mathbf{X}}$ is Hermitian, and so is $\log\cf{\mathbf{X}}$.

\emph{2)} We have:
\begin{align*}
 \cep{\mathbf{X}}{0} &= \frac{1}{2\pi} \int_{-\pi}^\pi \log \cf{\mathbf{X}} d\mu
 \stackrel{(a)}{=} \frac{1}{2\pi} \int_{-\pi}^\pi \ln |\cf{\mathbf{X}}| d\mu  \stackrel{(b)}{\le} 0
\end{align*}
where $(a)$ follows from~\eqref{eq:complex-log} and the fact that the phase of $\cf{\mathbf{X}}$ has odd symmetry; $(b)$ then follows from $|\cf{\mathbf{X}}|\le 1$.

\emph{3)} If the PMF is an even function, then $\cf{\mathbf{X}}$ and $\log\cf{\mathbf{X}}$ are real; the muculants have even symmetry by Fourier transform properties.

\emph{4)} If $\mathbb{E}(\mathbf{X})<\infty$, then $\cf{\mathbf{X}}$ is uniformly continuous and continuously differentiable; since $\log$ is piecewise continuous, we have that $\log\cf{\mathbf{X}}$ is piecewise continuous and piecewise differentiable; indeed,
\begin{equation}
 \frac{d}{d \mu} \log \cf{\mathbf{X}} = \frac{\cf{\mathbf{X}}'}{\cf{\mathbf{X}}}
\end{equation}
where $\cf{\mathbf{X}}$ vanishes on an at most countable set (otherwise the muculants would not exist, cf.~\eqref{eq:paleyWiener4}). Since the Fourier transform of a periodic, piecewise continuous and piecewise differentiable function converges pointwise~\cite[p.~105] {maths-sp11}, pointwise convergence of~\eqref{eq:complex-cepstrum-sum} follows. That, in this case, the muculants sum to zero, follows from evaluating~\eqref{eq:complex-cepstrum-sum} at $\mu=0$ and the fact that $\cf{\mathbf{X}}|_{\mu=0}=1$.

\emph{5)} That $\lim_{n\to \pm\infty} \cep{\mathbf{X}}{n}=0$ is a direct consequence of the Lebesgue-Riemann theorem and the fact that $\log\cf{\mathbf{X}}$ is absolutely integrable~\cite[p.~94]{maths-sp11}. If $\mathbb{E}(\mathbf{X})<\infty$, note that pointwise convergence implies bounded variation~\cite[p.70]{pereyra2012harmonic}; the result about the order of convergence follows from~\cite{taibleson}.

\emph{6)}
If $\mathbf{X}$ and $\mathbf{Y}$ are independent RVs, then $\cf{\mathbf{X}+\mathbf{Y}}=\cf{\mathbf{X}}\cdot \cf{\mathbf{Y}}$. The desired result follows from~\eqref{eq:logproperties} and the linearity of the Fourier series expansion.

\section*{Acknowledgment}
The authors would like to thank 
Franz Lehner, Institute of Mathematical Structure Theory, Graz University of Technology, 
Gennadiy Chystyakov, Faculty of Mathematics, Universit\"at Bielefeld, and
Paul Meissner  for fruitful discussions during the preparation of this manuscript.

\ifCLASSOPTIONcaptionsoff
  \newpage
\fi

 \bibliographystyle{IEEEtran}
 \bibliography{Muculants}

\begin{thebibliography}{10}
\providecommand{\url}[1]{#1}
\csname url@samestyle\endcsname
\providecommand{\newblock}{\relax}
\providecommand{\bibinfo}[2]{#2}
\providecommand{\BIBentrySTDinterwordspacing}{\spaceskip=0pt\relax}
\providecommand{\BIBentryALTinterwordstretchfactor}{4}
\providecommand{\BIBentryALTinterwordspacing}{\spaceskip=\fontdimen2\font plus
\BIBentryALTinterwordstretchfactor\fontdimen3\font minus
  \fontdimen4\font\relax}
\providecommand{\BIBforeignlanguage}[2]{{%
\expandafter\ifx\csname l@#1\endcsname\relax
\typeout{** WARNING: IEEEtran.bst: No hyphenation pattern has been}%
\typeout{** loaded for the language `#1'. Using the pattern for}%
\typeout{** the default language instead.}%
\else
\language=\csname l@#1\endcsname
\fi
#2}}
\providecommand{\BIBdecl}{\relax}
\BIBdecl

\bibitem{mendelHOS}
J.~M. Mendel, ``Tutorial on higher-order statistics (spectra) in signal
  processing and system theory: theoretical results and some applications.''
  \emph{Proceedings of the IEEE}, vol.~79, no.~3, pp. 278--305, 1991.

\bibitem{petropulu}
C.~L. Nikias and A.~P. Petropulu, \emph{Higher-Order Spectra Analysis: A
  Nonlinear Signal Processing Framework}, Englewood Cliffs, NJ, 1993.

\bibitem{moments-cumulants}
R.~A. Fisher and E.~A. Cornish, ``Moments and cumulants in the specification of
  distributions,'' \emph{Revue de l'Institut international de Statistique}, pp.
  307--320, 1938.

\bibitem{qMoments}
C.~Tsallis, A.~R. Plastino, and R.~F. Alvarez-Estrada, ``Escort mean values and
  the characterization of power-law-decaying probability densities.''
  \emph{Journal of Mathematical Physics}, vol.~50, no.~4, 2009.

\bibitem{qCumulants}
A.~Rodriguez and C.~Tsallis, ``A generalization of the cumulant expansion.
  application to a scale-invariant probabilistic model.'' \emph{Journal of
  Mathematical Physics}, vol.~51, no.~7, 2010.

\bibitem{cramer1957}
H.~Cramer, \emph{Mathematical Methods of Statistics}.\hskip 1em plus 0.5em
  minus 0.4em\relax Princeton University Press, 1957.

\bibitem{petrov1972}
V.~V. Petrov, \emph{Sums of Independent Random Variables}.\hskip 1em plus 0.5em
  minus 0.4em\relax Springer Verlag, 1975.

\bibitem{edgeworthExpansion}
S.~Blinnikov and R.~Moessner, ``Expansions for nearly {Gaussian}
  distributions,'' \emph{Astronomy and Astrophysics Supplement Series}, vol.
  130, pp. 193--205, 1998.

\bibitem{char_func_luk}
E.~Lukacs, \emph{Characteristic Functions. 2nd rev}.\hskip 1em plus 0.5em minus
  0.4em\relax Griffin, 1970.

\bibitem{Bogert_1963}
B.~P. Bogert, M.~J.~R. Healy, and J.~W. Tukey, ``The quefrency analysis of time
  series for echoes: cepstrum, pseudo-autocovariance, cross-cepstrum, and saphe
  cracking,'' in \emph{Proceedings Symp. Time Series Analysis}, M.~Rosenblatt,
  Ed.\hskip 1em plus 0.5em minus 0.4em\relax Wiley, 1963, pp. 209--243.

\bibitem{oppenheim-phd}
A.~V. Oppenheim, ``Superposition in a class of nonlinear systems.'' Ph.D.
  dissertation, Massachusetts Inst of Tech Camridge Research Lab of
  Electronics, 1965.

\bibitem{oppenheim_quefrency}
A.~V. Oppenheim and R.~Schafer, ``From frequency to quefrency: A history of the
  cepstrum,'' \emph{IEEE Signal Processing Magazine}, vol.~21, no.~5, pp.
  95--106, 2004.

\bibitem{schroeder-cepstrum-predictor}
M.~R. Schroeder, ``Direct (nonrecursive) relations between cepstrum and
  predictor coefficients,'' \emph{IEEE Transactions on Acoustics, Speech and
  Signal Processing}, vol.~29, pp. 297--301, 1981.

\bibitem{khare-moments}
A.~Khare and T.~Yoshikawa, ``Moment of cepstrum and its applications,''
  \emph{IEEE Transactions on Signal Processing,}, vol.~40, no.~11, pp.
  2692--2702, 1992.

\bibitem{luk-survey-cf}
E.~Lukacs, ``A survey of the theory of characteristic functions,''
  \emph{Advances in Applied Probability}, vol.~4, pp. 1--38, 1972.

\bibitem{applications-luk}
E.~Lukacs and R.~G. Laha, \emph{Applications of Characteristic
  Functions}.\hskip 1em plus 0.5em minus 0.4em\relax Griffin, 1964.

\bibitem{bisgaard2000}
T.~M. Bisgaard and Z.~Sasv\'ari, \emph{Characteristic Functions and Moment
  Sequences: Positive Definiteness in Probability}.\hskip 1em plus 0.5em minus
  0.4em\relax Nova Publishers, 2000.

\bibitem{tribolet-1977}
J.~M. Tribolet, ``A new phase unwrapping algorithm,'' \emph{IEEE Transactions
  on Acoustics, Speech and Signal Processing}, vol.~2, pp. 170--177, 1977.

\bibitem{brown-complex}
J.~W. Brown and R.~V. Churchill, \emph{Complex Variables and
  Applications}.\hskip 1em plus 0.5em minus 0.4em\relax McGraw Hill, 2009.

\bibitem{oppenheim2010discrete}
A.~V. Oppenheim and R.~W. Schafer, \emph{Discrete-Time Signal Processing},
  3rd~ed.\hskip 1em plus 0.5em minus 0.4em\relax Pearson Education, Limited,
  2010.

\bibitem{mc-gowan}
R.~McGowan and R.~Kuc, ``A direct relation between a signal time series and its
  unwrapped phase.'' \emph{IEEE Transactions on Acoustics, Speech and Signal
  Processing}, vol.~30, no.~5, pp. 719--726, 1982.

\bibitem{karam}
Z.~N. Karam, ``Computation of the one-dimensional unwrapped phase,'' Master's
  thesis, Massachusetts Institute of Technology, 2006.

\bibitem{mattner-cumulants}
L.~Mattner, ``What are cumulants?'' \emph{Documenta Mathematica}, vol.~4, pp.
  601--622, 1999.

\bibitem{rota-shen}
G.~C. Rota and J.~Shen, ``On the combinatorics of cumulants.'' \emph{Journal of
  Combinatorial Theory}, vol. Series A 91.1, pp. 283--304, 2000.

\bibitem{muculants-arxiv}
C.~Knoll, B.~C. Geiger, and G.~Kubin, ``The muculants--a higher-order
  statistical approach,'' \emph{arXiv preprint arXiv:1506.04518v1}, 2015.

\bibitem{papoulis_fourier}
A.~Papoulis, \emph{The Fourier Integral and its Applications}.\hskip 1em plus
  0.5em minus 0.4em\relax McGraw Hill, 1962.

\bibitem{papoulis_probability}
A.~Papoulis and S.~U. Pillai, \emph{Probability, Random Variables, and
  Stochastic Processes}, ser. McGraw-Hill electrical and electronic engineering
  series.\hskip 1em plus 0.5em minus 0.4em\relax McGraw-Hill, 2002, vol. 4th
  edition.

\bibitem{maths-sp11}
S.~B. Damelin and W.~{Miller, Jr.}, \emph{The Mathematics of Signal
  Processing}.\hskip 1em plus 0.5em minus 0.4em\relax Cambridge University
  Press, 2011.

\bibitem{good}
I.~J. Good, ``A new formula for cumulants,'' \emph{Mathematical Proceedings of
  the Cambridge Philosophical Society}, vol.~78, pp. 333--337, 1975.

\bibitem{negentropy}
C.~S. Withers and S.~Nadarajah, ``Negentropy as a function of cumulants,''
  Tech. Rep., 2011, research Report No. 15, Probability and Statistics Group,
  School of Mathematics, The University of Manchester.

\bibitem{widrow_quantization}
B.~Widrow, I.~Kollar, and M.~C. Liu, ``Statistical theory of quantization,''
  \emph{IEEE Transactions on Instrumentation and Measurement}, vol.~45, pp.
  353--361, 1996.

\bibitem{pereyra2012harmonic}
M.~C. Pereyra and L.~A. Ward, \emph{Harmonic analysis: from Fourier to
  wavelets}.\hskip 1em plus 0.5em minus 0.4em\relax American Mathematical Soc.,
  2012, vol.~63.

\bibitem{taibleson}
M.~Taibleson, ``Fourier coefficients of functions of bounded variation,''
  \emph{Proceedings of the American Mathematical Society}, vol.~18, no.~4, p.
  766, 1967.

\end{thebibliography}
%

%
%
%




\end{document}